\magnification 1200
\input vladimi
\hsize 17truecm
\vsize 24truecm
\hoffset 0truecm
\voffset -0.5truecm
\parindent 1em
\tolerance 9000
\tenpoint\frenchspacing
\leftline{\ }\vskip 0.25cm
{\leftskip 0cm plus 1fill\rightskip 0cm plus 1fill\parindent 0cm\baselineskip 15pt
\TITLE Осцилляционные свойства одной многоточечной граничной задачи
четвёртого порядка со спектральным параметром в граничном условии
\par\vskip 0.15cm\rm А.$\,$А.~Владимиров, Е.$\,$С.~Карулина\par}
\vskip 0.25cm
$$
	\vbox{\hsize 0.75\hsize\leftskip 0cm\rightskip 0cm
	\eightpoint\rm
	{\bf Аннотация:\/} Для одной многоточечной граничной задачи четвёртого
	порядка со спектральным параметром в граничном условии установлена
	простота спектра и наличие осцилляционных свойств у системы производных
	собственных функций.\par
	}
$$

\vskip 0.25cm
\section{Введение}
\subsection
Рассмотрим граничную задачу
$$
	\displaylines{(py'')''-(qy')'=\lambda ry,\cr
	y(0)=y'(0)=0,\cr
	y(\xi_i+0)-y(\xi_i-0)=y'(\xi_i+0)-\eta_iy'(\xi_i-0)={}\kern 5cm\cr
	{}=\eta_i(py'')(\xi_i+0)-(py'')(\xi_i-0)-\alpha_iy'(\xi_i-0)={}\cr
	\kern 5cm{}=[(py'')'-qy'](\xi_i+0)-[(py'')'-qy'](\xi_i-0)=0,\cr
	(py'')(1)-\alpha\lambda y'(1)=[(py'')'-qy'](1)+\beta\lambda y(1)=0,}
$$
где $\lambda$~--- спектральный параметр, $0<\xi_0<\ldots<\xi_m<1$, $\eta_i>0$,
$\alpha\geqslant 0$, $\beta>0$, $p,1/p\in L_\infty[0,1]$ неотрицательны,
$q\in L_1[0,1]$ вещественнозначна, а $r\in L_1[0,1]$ почти всюду положительна.
Основным результатом настоящей статьи является следующее утверждение:

\subsubsection\label{prop:1.1}
{\it В случае, если все собственные значения
$$
	\lambda_0\leqslant\lambda_1\leqslant\ldots\leqslant\lambda_n\leqslant\ldots
$$
рассматриваемой задачи положительны, они являются простыми, причём
для соответствующих собственных функций $y_n$ производные $y'_n$ имеют в точности
по $n$ перемен знака.}

\medskip
В некоторых частных случаях аналогичные утверждения были установлены, например,
в работе [\cite{KA:2006}].

\subsection
Операторную модель, на основе которой будет проводиться изучение поставленной
задачи, мы вводим следующим образом. Пусть $\tau\colon[0,1]\to[0,1]$~--- непрерывное
возрастающее кусочно-линейное отображение отрезка $[0,1]$ на себя, имеющее
изломы только в точках $\xi_i$ и подчиняющееся при этом уравнениям
$$
	\tau'(\xi_i+0)=\eta_i\tau'(\xi_i-0).
$$
Тогда всякая функция $y$, удовлетворяющая граничным условиям, допускает
представле\-ние в виде $y=u\circ\tau$, где
$$
	u\in\frak D\rightleftharpoons\{u\in W_2^2[0,1]\;:\; u(0)=u'(0)=0\}.
$$
При этом решению задачи отвечает такая и только такая функция $u\in\frak D$,
которая при всяком выборе пробной функции $v\in\frak D$ подчиняется равенству
$$
	\int_0^1 \hat pu''\overline{v''}\,dx+\langle\hat q,\overline{u'}v'\rangle-
		\lambda\cdot\left[\int_0^1\hat ru\overline{v}\,dx+
		\hat\alpha u'(1)\overline{v'(1)}+\beta u(1)\overline{v(1)}\right]=0,
$$
где $\hat p\circ\tau=p\cdot[\tau']^3$, $\hat r\circ\tau=r/\tau'$, $\hat\alpha=
\alpha\cdot[\tau'(1)]^2$, а обобщённая функция $\hat q\in W_2^{-1}[0,1]$
определяется тождеством
$$
	\langle\hat q,u\rangle\equiv\int_0^1 q\cdot[\tau']^2\,
		\overline{u\circ\tau}\,dx+\sum_{i=0}^m\alpha_i\cdot
		[\tau'(\xi_i-0)]^2\,\overline{u(\tau(\xi_i))}.
$$
Соответственно, исходная спектральная граничная задача равносильна спектральной
задаче для линейного операторного пучка $\hat T\colon\Bbb C\to{\cal B}(\frak D,
\frak D^*)$ вида
$$
	\langle\hat T(\lambda)u,v\rangle\equiv\int_0^1 \hat pu''\overline{v''}\,dx+
		\langle\hat q,\overline{u'}v'\rangle-
		\lambda\cdot\left[\int_0^1\hat ru\overline{v}\,dx+
		\hat\alpha u'(1)\overline{v'(1)}+\beta u(1)\overline{v(1)}\right].
$$
Очевидное совпадение знаков величин $y'(x)=u'(\tau(x))\cdot\tau'(x)$
и $u'(\tau(x))$ в каждой точке $x\in [0,1]$, не являющейся точкой излома функции
$\tau$, означает, что искомые осцилляционные свойства собственных функций исходной
задачи также совпадают с таковыми для собственных функций пучка $\hat T$.
Именно этот пучок и будет подвергаться изучению в основной части статьи.
При этом выполнение предположений утверждения~\ref{prop:1.1} будет предполагаться
без специальных оговорок.

Отметим, что в работе [\cite{Ku:2017}] многоточечные граничные задачи
рассмотренного нами вида названы «задачами на графах-цепочках». Построенная
в настоящем пункте операторная модель показывает, что с принципиальной точки
зрения такие задачи не отличаются от обычных двухточечных граничных задач.


\section{Редукция задачи и простота спектра}
\subsection
Заметим, что оператор $\hat T(0)$ представляет собой вполне непрерывное возмущение
некоторого равномерно положительного оператора, а оператор $\hat T'(0)$ есть
отрицательный вполне непрерывный оператор. Соответственно, из части~5 теоремы
[\cite{LSY:1993}: Theorem~1] немедленно вытекает, что оператор $\hat T(0)$
является положительно определённым.

Введём в рассмотрение пространство
$$
	\frak D_1\rightleftharpoons\left\{u\in W_2^1[0,1]\;:\; u(0)=0\right\},
	\leqno(\equation)
$$\label{eq:4}%
а также биекцию $Q\colon\frak D_1\to\frak D$ вида $[Qy]'\equiv y$. Оператор
$S\rightleftharpoons Q^*\hat T(0)Q$ при этом положительно определён и подчиняется
тождеству
$$
	\langle Su,v\rangle\equiv\int_0^1\hat pu'\overline{v'}\,dx+
		\langle\hat q,\overline{u}v\rangle.
$$
Ввиду известного (см., например, [\cite{Vl:2016}:~\S$\,$3]) факта знакорегулярности
положительно определённых операторов Штурма--Лиувилля, это означает существование
равномерно положительной функции $\sigma\in W_2^1[0,1]$ со свойством
$$
	(\forall v\in\frak D_1)\qquad \int_0^1\hat p\sigma'\overline{v'}\,dx
		+\langle\hat q,\sigma v\rangle=\gamma\overline{v(1)},
	\leqno(\equation)
$$\label{eq:1}%
где $\gamma>0$~--- некоторая постоянная. Полагая без ограничения общности
выполненным равенство $\int_0^1\sigma\,dx=1$, мы можем ввести в рассмотрение
биективное отображение $\omega\colon[0,1]\to [0,1]$ со свойством $\omega'=\sigma$,
а также связанное с ним биективное отображение $V:\frak D\to\frak D$ вида 
$$
	\leqalignno{[Vu]&=u\circ\omega,\cr
		[Vu]'&=(u'\circ\omega)\cdot\sigma,\cr 
		[Vu]''&=(u''\circ\omega)\cdot\sigma^2+(u'\circ\omega)\cdot\sigma'.}
$$
Количества знакоперемен функций $u'$ и $[Vu]'$ при этом очевидным образом совпадают.
Соответственно, для доказательства утверждения~\ref{prop:1.1} достаточно установить
искомое для случая пучка операторов $\tilde T(\lambda)\rightleftharpoons
V^*\hat T(\lambda)V$, имеющих вид
$$
	\langle\tilde T(\lambda)u,v\rangle\equiv 
		\int_0^1\tilde pu''\overline{v''}\,dx+
		\gamma\sigma(1) u'(1)\overline{v'(1)}-\lambda\cdot\left[\int_0^1\tilde r
		u\overline{v}\,dx+\tilde\alpha u'(1)\overline{v'(1)}+
		\beta u(1)\overline{v(1)}\right],
$$
где $\tilde p\circ\omega=\hat p\sigma^3$, $\tilde r\circ\omega=\hat r/\sigma$
и $\tilde\alpha=\hat\alpha\sigma^2(1)$, а $\gamma$~--- постоянная
из тождества~\eqref{eq:1}. Собственные значения и собственные функции такого
пучка определяются классически понимаемой граничной задачей
$$
	\belowdisplayskip 2pt
	(\tilde pu'')''=\lambda\tilde ru,\leqno(\equation)
$$\label{eq:2}%
$$
	\abovedisplayskip 2pt
	u(0)=u'(0)=(\tilde pu'')(1)+[\gamma\sigma(1)-\tilde\alpha\lambda]u'(1)=
		(\tilde pu'')'(1)+\beta\lambda u(1)=0.\leqno(\equation)
$$\label{eq:3}%

Отметим, что проведённый в настоящем пункте переход от пучка $\hat T$ к пучку
$\tilde T$ ранее применялся, например, в работах [\cite{Vl:2006}, \cite{BVS:2013}].

\subsection
Имеет место следующий факт:

\subsubsection\label{prop:2.2}
{\it Всякое собственное значение граничной задачи~\eqref{eq:2}, \eqref{eq:3}
является простым. При этом всякая собственная функция указанной задачи
подчиняется неравенству $(\tilde pu'')(x)\neq 0$ при $x=0$, а также
при $u'(x)=0$.
}
\proof
Ввиду равносильности рассматриваемой граничной задачи спектральной задаче для пучка
$\tilde T$, собственные значения этой задачи положительны. Соответственно,
собственные функции подчинены известной лемме [\cite{LN:1958}: Lemma~2.1].
А именно, любое нетривиальное решение граничной задачи~\eqref{eq:2}, \eqref{eq:3},
удовлетворяющее дополнительному условию $(\tilde pu'')(0)=0$, обязано иметь
отличные от нуля значения величин $u(1)$ и $(\tilde pu'')'(1)$, знак которых
совпадает со знаком величины $(\tilde pu'')'(0)$. Однако это несовместимо
с условием $\beta>0$. Таким образом, нами установлен факт заведомой простоты
собственных значений, а также выполнение неравенства $(\tilde pu'')(0)\neq 0$
для всякой собственной функции.

Далее, каждая точка $x\in (0,1)$ со свойством $u'(x)=(\tilde pu'')(x)=0$ должна
подчиняться какому-то из неравенств $(\tilde pu'')'(x)\cdot u(x)\geqslant 0$ либо
$(\tilde pu'')'(x)\cdot u(x)\leqslant 0$. В первом случае вышеуказанная лемма снова
гарантирует выполнение несовместимого с условием $\beta>0$ неравенства
$(\tilde pu'')'(1)\cdot u(1)>0$. Во втором же случае из леммы [\cite{LN:1958}:
Lemma~2.2] вытекает противоречащее граничным условиям неравенство $u(0)\neq 0$.
\endproof


\section{Завершение доказательства}
\subsection
Введём в рассмотрение действующий в пространстве $\frak D_1$ вида~\eqref{eq:4}
оператор $K\rightleftharpoons -Q^{-1}[\tilde T(0)]^{-1}\tilde T'(0)Q$, где $Q\colon
\frak D_1\to\frak D$ есть, как и ранее, биективный оператор интегрирования.
Имеет место следующий факт:

\subsubsection\label{prop:3.1}
{\it Оператор $K$ не повышает числа знакоперемен никакой вещественной функции
$u\in\frak D_1$.
}
\proof
Равенство $Ku=v$ равносильно равенству $-\tilde T'(0)Qu=\tilde T(0)Qv$, а потому
и равенству $-Q^*\tilde T'(0)Qu=Q^*\tilde T(0)Qv$. Оператор $Q^*\tilde T(0)Q\colon
\frak D_1\to\frak D_1^*$ при этом представляет собой положительно определённый
оператор Штурма--Лиувилля, и потому знакорегулярен (см., например,
[\cite{Vl:2016}]). Соответственно, число знакоперемен функции $v\in\frak D_1$
мажорируется числом знакоперемен обобщённой функции $Q^*\tilde T(0)Qv\in
\frak D_1^*$. Для завершения доказательства, таким образом, остаётся установить,
что число знакоперемен обобщённой функции $-Q^*\tilde T'(0)Qu\in\frak D_1^*$
мажорируется числом знакоперемен функции $u\in\frak D_1$.

Заметим, что имеет место равенство
$$
	-Q^*\tilde T'(0)Qu=w+\tilde\alpha u(1)\bolddelta_1,\leqno(\equation)
$$\label{eq:5}%
где $\bolddelta_1\in\frak D_1^*$~--- дельта-функция с носителем в точке $1$,
а функция $w\in C[0,1]$ определяется тождеством
$$
	w(t)\equiv\int_t^1\tilde r\cdot(Qu)\,dx+\beta\cdot(Qu)(1).
$$
Непосредственно из интегральной теоремы о среднем следует, что число знакоперемен
функции $w$ не превосходит такового для функции $Qu$, а тогда и для функции $u$.
Соответственно, обобщённая функция~\eqref{eq:5} могла бы иметь большее,
чем функция $u$, число знакоперемен лишь в случае $u(1)\neq 0$ и $(Qu)(1)\cdot u(1)
\leqslant 0$. Но в этом случае число знакоперемен функции $Qu$ заведомо строго
меньше, чем такое число для функции $u$. Тем самым, число знакоперемен обобщённой
функции~\eqref{eq:5} мажорируется числом знакоперемен функции $u$ во всех
возможных случаях.
\endproof

\subsection
Легко видеть, что всякой собственной паре $\{\lambda,u\}$ операторного пучка
$\tilde T$ отвечает собственная пара $\{\lambda^{-1},Q^{-1}u\}$ оператора $K$,
и наоборот. Соответственно, для завершения доказательства утверждения~\ref{prop:1.1}
достаточно установить, что система $\{h_n\}_{n=0}^\infty$ собственных функций
оператора $K$, занумерованных в порядке убывания соответствующих собственных
значений, обладает стандартными осцилляционными свойствами.

Заметим, что система $\{Qh_n\}_{n=0}^\infty$ собственных функций пучка $\tilde T$
образует базис Рисса в пространстве $\frak D$. Ввиду биективности оператора
$Q$ это означает, что система $\{h_n\}_{n=0}^\infty$ образует базис Рисса
в пространстве $\frak D_1$. Такое наблюдение, в свою очередь, позволяет провести
следующие два рассуждения.

Во-первых, при любом выборе числа $n\in\Bbb N$ внутри линейной оболочки набора
$\{x^m\}_{m=1}^{n+1}$ функций класса $\frak D_1$ существует многочлен вида
$g=\sum\limits_{k=N}^{\infty}{c_kh_k}$, где $N\geqslant n$ и $c_N\neq 0$. При этом
будет справедливо равенство $\lim_{m\to\infty}\lambda_N^mK^mg=c_Nh_N$, ввиду
утверждения~\ref{prop:3.1} означающее, что число знакоперемен функции $h_N$
не превосходит такового для многочлена $g$. Иначе говоря, для любого значения
$n\in\Bbb N$ существует значение $N\geqslant n$, для которого число знакоперемен
собственной функции $h_N$ не будет превосходить $n$.

Во-вторых, при любом выборе числа $n\in\Bbb N$ внутри линейной оболочки набора
$\{h_m\}_{m=0}^n$ найдётся функция $g=\sum\limits_{k=0}^n {c_kh_k}$,
удовлетворяющая неравенству $c_n\neq 0$ и имеющая не менее $n$~перемен знака.
При этом каждая из функций вида $g_m=\sum\limits_{k=0}^n(\lambda_k/\lambda_n)^m
c_kh_k$, обладающих свойством $\lambda_n^m K^mg_m=g$, также должна, согласно
утверждению~\ref{prop:3.1}, иметь не менее $n$ перемен знака. Между тем,
функциональные последовательности $\{g_m\}_{m=0}^\infty$
и $\{\tilde pg_m'\}_{m=0}^\infty$ равномерно сходятся к функциям $c_nh_n$
и $c_n\tilde ph_n'$, соответственно. Согласно утверждению~\ref{prop:2.2}
это означает, что при достаточно больших значениях индекса $m\in\Bbb N$ функции
$g_m$ имеют совпадающее с таковым для функции $h_n$ число знакоперемен.
Иначе говоря, для любого значения $n\in\Bbb N$ число знакоперемен собственной
функции $h_n$ не может быть меньшим, чем $n$.

Объединяя результаты проведённых двух рассуждений, получаем, что при любом
выборе значения $n\in\Bbb N$ соответствующая функция $h_n$ имеет в точности
$n$ знакоперемен. Доказательство утверждения~\ref{prop:1.1} тем самым завершено.

\subsection
Утверждение~\ref{prop:1.1} может быть легко дополнено утверждениями о поведении
числа знакоперемен прочих квазипроизводных собственных функций изучаемой задачи.
Например, имеет место следующий факт:

\subsubsection
{\it При любом выборе индекса $n\in\Bbb N$ соответствующая собственная функция
$y_n$ имеет не менее $n-1$ и не более $n$~перемен знака. При этом может быть
указано такое не зависящее от выбора параметра $\alpha\geqslant 0$ число
$\varkappa>0$, что для всякого собственного значения $\lambda_n$ со свойством
$\alpha\lambda_n\leqslant\varkappa$ соответствующая собственная функция $y_n$
будет иметь в точности $n$~перемен знака.
}
\proof
Как и ранее, достаточно рассмотреть случай граничной задачи~\eqref{eq:2},
\eqref{eq:3}. Ввиду того, что функция $u_n'$ имеет в точности $n$~знакоперемен,
функция $u_n$ заведомо не может иметь более $n$~таковых. С другой стороны,
функция $\tilde pu_n''$ имеет не менее $n$~знакоперемен, а функция
$(\tilde pu_n'')'$~--- не менее $n-1$. Потому из уравнения~\eqref{eq:2} и заданной
граничными условиями~\eqref{eq:3} связи между значениями $(\tilde pu_n'')'(1)$
и $u_n(1)$ вытекает, что функция $u_n$ имеет не менее $n-1$~знакоперемен.

Предположим теперь выполненным неравенство $\tilde\alpha\lambda_n\leqslant
\gamma\sigma(1)$. Здесь в случае $(\tilde pu_n'')(1)\neq 0$ функция $\tilde pu_n''$
имеет не менее $n+1$~знакоперемен. Это означает заведомое наличие не менее
$n$~знакоперемен у функции $(\tilde pu_n'')'$, а тогда и у функции $u_n$.
\endproof

\vskip 0.4cm
\eightpoint\rm
{\leftskip 0cm\rightskip 0cm plus 1fill\parindent 0cm
\bf Литература\par\penalty 20000}\vskip 0.4cm\penalty 20000
\bibitem{KA:2006} Керимов~Н.$\,$Б., Алиев~З.$\,$С. {\it Базисные свойства
одной спектральной задачи со спектральным параметром в граничном условии}~//
Матем.~сб.~--- 2006.~--- Т.~197, \No~10.~--- С.~65--86.
\bibitem{Ku:2017} Кулаев~Р.$\,$Ч. {\it К вопросу о неосцилляции дифференциального
уравнения на графе}~// Владикавказ. матем.~журнал.~--- 2017.~--- Т.~19, \No~3.~---
С.~31--40.
\bibitem{LSY:1993} Lancaster~P., Shkalikov~A., Qiang~Ye. {\it Strongly
definitizable linear pencils in Hilbert space}~// Integr. Equat. Oper. Th.~---
1993.~--- V.~17.~--- P.~338--360.
\bibitem{Vl:2016} Владимиров~А.$\,$А. {\it К вопросу об осцилляционных свойствах
положительных дифференциальных операторов с сингулярными коэффициентами}~//
Матем. заметки.~--- 2016.~--- Т.~100, \No~6.~--- С.~800--806.
\bibitem{Vl:2006} Бен~Амара~Ж., Владимиров~А.$\,$А. {\it Об осцилляции собственных
функций задачи четвёртого порядка со спектральным параметром в граничном условии}~//
Фунд. и прикл. матем.~--- 2006.~--- Т.~12, \No~4.~--- С.~41--52.
\bibitem{BVS:2013} Бен~Амара~Ж., Владимиров~А.$\,$А., Шкаликов~А.$\,$А.
{\it Спектральные и осцилляционные свойства одного линейного пучка дифференциальных
операторов четвёртого порядка}~// Матем.~заметки.~--- 2013.~--- Т.~94, \No~1.~---
С.~55--67.
\bibitem{LN:1958} Leighton~W., Nehari~Z. {\it On the oscillation of solutions
of self-adjoint linear differential equations of the fourth order}~//
Trans. Amer. Math. Soc. --- 1958. --- Vol. 89. --- P.~325--377.
\bye